\documentclass[11pt,a4paper]{article}
\usepackage[T1]{fontenc}
\usepackage[utf8]{inputenc}
\usepackage[english]{babel}
\usepackage{lmodern,microtype}
\usepackage{amsmath,amssymb,amsthm}
\usepackage{geometry,xcolor,hyperref}
\hypersetup{colorlinks=true,linkcolor=blue!45!black,citecolor=blue!45!black,
urlcolor=blue!45!black,pdftitle={Online packing of the three largest trees}}
\newtheorem{theorem}{Theorem}[section]
\newtheorem{lemma}[theorem]{Lemma}
\theoremstyle{remark}

\title{\textbf{Online packing of the three largest trees}}
\author{%
Maciej Cisiński\thanks{The author was partially supported
by the Polish Ministry of Science and Higher Education.}\\
{\small Department of Discrete Mathematics, AGH University of Krakow}\\
{\small al. Mickiewicza 30,
30-059 Kraków, Poland}\\
{\small \texttt{cisinski@agh.edu.pl}}
}
\begin{document}
\maketitle
\begin{abstract}
For every $n\ge4$, arbitrary trees of orders $n,n-1,n-2$ can be
packed into $K_n$ online, in decreasing order of size: each tree
must be embedded before the next is revealed, and earlier embeddings
cannot be changed. This strengthens the classical three-tree packing
theorem of Hobbs, Bourgeois and Kasiraj (1987). In particular, the
second tree is embedded without knowing the third, and the resulting
pair of copies can be extended to accommodate every possible third
tree.

Our strategy is to embed the first two trees so that their union is
$2$-degenerate. We prove that the complement of every $2$-degenerate
graph on $n$ vertices contains every tree of order $n-2$, which
guarantees the final move. The main ingredient is an embedding lemma
for spiders whose center must occupy a prescribed vertex having
degree at most two in the forbidden graph. We then extend this
lemma to arbitrary trees by induction.
\end{abstract}


\section{Introduction}

All graphs are finite and simple. The order of a graph is its number
of vertices, and embeddings need not be induced. A family of graphs
\emph{packs} into $K_n$ if $K_n$ contains pairwise edge-disjoint
copies of its members; these copies may share vertices.

The classical tree packing conjecture of Gy\'arf\'as, studied by
Gy\'arf\'as and Lehel~\cite{GyarfasLehel1978}, asserts that any family
of trees $T_1,\ldots,T_n$ with $|V(T_i)|=i$ packs into $K_n$.
Since these trees have
\[
 \sum_{i=1}^{n}(i-1)=\binom n2
\]
edges in total, such a packing would decompose the edges of $K_n$.

Even packing only the largest trees in such a family is a substantial
problem. Bollob\'as~\cite{Bollobas1995} proposed the following
weakening of the tree packing conjecture: for every fixed $k$,
there exists an $n_0(k)$ such that, whenever $n\ge n_0(k)$,
arbitrary trees of orders $n,n-1,\ldots,n-k+1$ pack into $K_n$.
Hobbs, Bourgeois and Kasiraj~\cite{HobbsBourgeoisKasiraj1987}
proved that trees of orders $n,n-1,n-2$ always pack into $K_n$.
For sufficiently large $n$, \.{Z}ak~\cite{Zak2017} established the
corresponding result for the five largest trees.
Cisi\'nski and \.{Z}ak~\cite{CisinskiZak2024} extended this to the
seven largest trees, again for sufficiently large $n$.
Janzer and Montgomery~\cite{JanzerMontgomery2026} proved
Bollob\'as's conjecture in a stronger form, showing that a linear
number of the largest trees can be packed.

These are \emph{offline} packing problems: all trees are known before
any embedding is chosen. Here we consider an \emph{online} version
for the three largest trees. The host graph $K_n$ is fixed and known
to both the player and the adversary. The game proceeds as follows.
The adversary first reveals a tree $T_n$ of order $n$, and the player
embeds it into $K_n$. The adversary then reveals a tree $T_{n-1}$ of
order $n-1$, which the player must embed using only edges not used
by the first copy. Finally, the adversary reveals a tree $T_{n-2}$
of order $n-2$, and the player must embed it using the remaining
edges. Each embedding is fixed before the next tree is revealed
and cannot subsequently be changed. The adversary may choose each
new tree after observing all earlier embeddings. Thus the online
condition concerns the arrival of whole trees, rather than individual
vertices or edges.

The distinction from the offline problem is a change in the order
of quantifiers. In an offline packing, the embeddings of the first
two trees may depend on the third tree. In the online problem, the
first embedding must be chosen without knowing either later tree.
Once the second tree is revealed, its embedding must be chosen
without knowing the third. The same pair of embedded copies must
therefore admit an extension for every possible choice of the third
tree. An offline existence theorem alone does not guarantee this.

Our main result shows that this additional requirement can always
be met.

\begin{theorem}\label{thm:main}
For every $n\ge4$, there is an online strategy that packs arbitrary
trees $T_n,T_{n-1},T_{n-2}$ of orders $n,n-1,n-2$, respectively,
into $K_n$, when they are revealed in this order. Each tree is
embedded before the next is revealed, and no earlier embedding
is changed.
\end{theorem}

The strategy is to control the structure of the edges used by the
first two trees. A graph is \emph{$2$-degenerate} if every nonempty
subgraph has a vertex of degree at most two. Equivalently, its
vertices can be removed one at a time so that each removed vertex
has at most two neighbors among the vertices still present.
In particular, every subgraph of a $2$-degenerate graph is itself
$2$-degenerate.

We will show that, after placing the first tree arbitrarily, the
second tree can be embedded so that the union $H$ of the two copies
is $2$-degenerate. The unused edges then form the complement
$\overline H$. The following structural embedding theorem guarantees
that these edges can accommodate any third tree.

\begin{theorem}\label{thm:universal}
For every $n\ge4$ and every $2$-degenerate graph $H$ on $n$ vertices,
the complement $\overline H$ contains a copy of every tree of order
$n-2$.
\end{theorem}

Theorem~\ref{thm:universal} applies to every $2$-degenerate forbidden
graph, whether or not it arises as the union of two trees.
\section{Preliminary lemmas}

To use Theorem~\ref{thm:universal} in our online strategy, we first show that the first two
trees can always be packed with a $2$-degenerate union, while keeping
the first copy fixed. The following lemma uses the leaf-removal construction from the proof of Theorem 1 in \cite{GyarfasLehel1978}, with the additional
observation that the resulting union is $2$-degenerate.

\begin{lemma}[Pair preparation]\label{lem:pair}
Let $T$ be any spanning tree already embedded in $K_n$, where $n\ge2$.
Every tree $S$ of order $n-1$ can be embedded edge-disjointly from $T$
so that the union of the two copies is $2$-degenerate.
\end{lemma}

\begin{proof}
We induct on $n$. The case $n=2$ is immediate.

\smallskip\noindent
\textbf{Case 1: at least one of the trees is a star.}
If $T$ is a star, embed $S$ on the vertices other than its center.
If $S$ is a star, place its center at a leaf $x$ of $T$, leave the
neighbor of $x$ unused by $S$, and place its leaves on the remaining
vertices. In either construction, the union consists of a tree on
$n-1$ vertices together with an additional vertex adjacent to all
of them. Successively removing leaves of that tree, and then the
additional vertex, gives a $2$-degeneracy order.

\smallskip\noindent
\textbf{Case 2: neither tree is a star.}
Each tree has two leaves with distinct neighbors: take the endpoints
of a longest path. Choose such leaves $x_1,x_2$ in $T$ and $y_1,y_2$
in $S$. Delete them and apply induction to the remaining trees,
of orders $n-2$ and $n-3$, on the host vertices outside
$\{x_1,x_2\}$.

Let $a_1\ne a_2$ be the neighbors of $x_1,x_2$ in $T$, and let
$b_1\ne b_2$ be the images of the neighbors of $y_1,y_2$ in $S$.
First try placing $y_1$ at $x_1$ and $y_2$ at $x_2$.
If neither restored edge coincides with an edge of $T$, we are done.
Otherwise, $a_1=b_1$ or $a_2=b_2$. In either case, the inequalities
$a_1\ne a_2$ and $b_1\ne b_2$ imply that $a_1\ne b_2$ and
$a_2\ne b_1$. We can therefore place $y_1$ at $x_2$ and $y_2$ at
$x_1$, avoiding both conflicts.

Each restored host vertex has degree two in the union, and the two
are nonadjacent. Removing them first and then following a degeneracy
order of the smaller union proves the claim.
\end{proof}
The embedding of the first tree remains fixed throughout the
construction, and the second embedding is chosen without any
knowledge of the third tree.

\begin{lemma}[A rooted local extension]\label{lem:seven}
Let $H$ be a $2$-degenerate graph on seven vertices, and let
$G=\overline H$. If $G$ contains a four-vertex path starting at $r$,
then it also contains a five-vertex path starting at $r$.
\end{lemma}

\begin{proof}
We first establish two elementary facts.

\smallskip\noindent
\textbf{Fact 1.}
In the complement of a forest, a path with at least two vertices
and at least two unused vertices can be replaced by a path with
one more vertex and the same starting vertex.

\begin{proof}[Proof of Fact 1]
\renewcommand{\qedsymbol}{}
Let $x,y$ be the last two vertices of the path, and let $u,v$
be unused vertices. If $y$ has a complement neighbor among $u,v$,
append it. Otherwise $yu,yv$ are forest edges, so $uv$ is a
complement edge. If $xu$ or $xv$ is a complement edge, replace
$y$ by $u,v$ in the appropriate order. If neither is, the forest
contains the cycle $xuyvx$, a contradiction.
\end{proof}

\smallskip\noindent
\textbf{Fact 2.}
Suppose disjoint sets $X,Y$ are completely joined in $H$.
If $|Y|=1$, then $H[X]$ is a forest; if $|Y|=2$, then $H[X]$
has at most one edge.

\begin{proof}[Proof of Fact 2]
\renewcommand{\qedsymbol}{}
If $H[X]$ contains a cycle and $|Y|=1$, consider the subgraph
induced by the cycle vertices together with $Y$. Each cycle
vertex has two cycle neighbors and a neighbor in $Y$, while
the vertex of $Y$ has at least three neighbors on the cycle.
This contradicts $2$-degeneracy.

Now let $|Y|=2$. If $H[X]$ contains two distinct edges, let $U$
be the set of their endpoints. Then $|U|\ge3$, and every vertex
of $U$ has a neighbor in $U$. In $H[U\cup Y]$, each vertex of $U$
has that neighbor and both vertices of $Y$ as neighbors.
Each vertex of $Y$ is adjacent to all of $U$. Thus this subgraph
has minimum degree at least three, again a contradiction.
\end{proof}

\smallskip
Suppose now that no five-vertex path starts at $r$.

\smallskip\noindent
\textbf{First, $G$ is connected.}
Let $C$ be the component containing $r$. The given path ensures
that $|C|\ge4$. If $|C|=4$, all edges between $C$ and the other
three vertices belong to $H$, giving a $K_{3,3}$.
If $|C|=5$, Fact 2 gives $e(H[C])\le1$, so $G[C]$ has
a Hamilton path starting at $r$. If $|C|=6$, Fact 2 says
that $H[C]$ is a forest, and Fact 1 lengthens the given
path inside $C$. Each possibility is a contradiction, so $|C|=7$.

\smallskip\noindent
\textbf{Next, we determine the shape of a depth-first search tree.}
Take a DFS spanning tree $D$ of $G$, rooted at $r$.
Every edge of $G$ joins an ancestor to a descendant in $D$.
Thus subtrees rooted at distinct children of the same vertex
have no $G$-edges between them and are completely joined in $H$.
Also, $D$ has height at most three, since a vertex at depth four
would give a five-vertex path from $r$.

Some child subtree $A$ of $r$ has at least four vertices.
Otherwise, in $H-r$ every vertex has at least three neighbors
in the other child subtrees, contradicting $2$-degeneracy.
Let $a$ be the root of $A$. If $|A|=4$, Fact 2 gives
$e(H[A])\le1$, so $G[A]$ has a four-vertex path starting at $a$.
If $|A|=5$, then $H[A]$ is a forest; start with an edge of $D$
from $a$ and apply Fact 1 twice to obtain such a path.
Prepending $r$ gives a contradiction in either case.

Hence $|A|=6$, and $a$ is the only child of $r$.
Among the five vertices below $a$, some child subtree $B$ of $a$
has at least three vertices. Otherwise their induced subgraph
in $H$ has minimum degree at least three.
The height bound makes $B$ a star with center $b$.
Its leaves are pairwise adjacent in $H$. If $|B|\ge5$, four leaves
give a $K_4$ in $H$; if $|B|=4$, its three leaves together with
a vertex in another child subtree of $a$ do so.
Thus $B$ has exactly three vertices, with leaves $c,d$.

The remaining two vertices must form one child subtree of $a$,
with tree edge $ef$ and $e$ adjacent to $a$.
Indeed, two separate children, together with $c,d$, would give
a $K_4$ in $H$. Consequently,
\[
 E(D)=\{ra,ab,bc,bd,ae,ef\}.
\]

\smallskip\noindent
\textbf{Finally, this shape also gives a contradiction.}
Each of $b,c,d,e,f$ starts a four-vertex path in $D-r$:
the corresponding vertex sequences are
\[
 (b,a,e,f),\quad (c,b,a,e),\quad (d,b,a,e),\quad
 (e,a,b,c),\quad (f,e,a,b).
\]
A $G$-edge from $r$ to any of these vertices would therefore
give a five-vertex path starting at $r$. Hence all five are
$H$-neighbors of $r$.

Moreover, $\{c,d\}$ is completely joined to $\{e,f\}$ in $H$,
and $cd\in E(H)$, by the DFS property.
Since $ef\in E(G)$, the vertices $\{r,c,d,e,f\}$ induce $K_5-ef$
in $H$. This graph has minimum degree three, contradicting
$2$-degeneracy.
\end{proof}

\begin{lemma}[Paths and terminal repair]\label{lem:paths}
Let $H$ be a $2$-degenerate graph on $n\ge4$ vertices, and let
$G=\overline H$. Then:
\begin{enumerate}
\item $G$ contains a path on $n-2$ vertices.
\item Suppose a tree $T$ contains a path $v_0v_1v_2v_3v_4$,
where $v_1,v_2,v_3$ have degree two in $T$ and $v_4$ is a leaf.
Every embedding of $T-v_4$ into $G$ leaving exactly three host
vertices unused can be modified to embed $T$, changing only
the images of $v_1,v_2,v_3$ and adding an image for $v_4$.
In particular, the image of $v_0$ remains fixed.
\end{enumerate}
\end{lemma}

\begin{proof}
We first prove the second assertion.
Consider the subgraph of $G$ induced by the images of
$v_0,v_1,v_2,v_3$ together with the three unused host vertices.
Its complement is an induced subgraph of $H$ and is therefore
$2$-degenerate. Lemma~\ref{lem:seven} gives a five-vertex path
starting at the image of $v_0$.
Use this path to embed $v_0v_1v_2v_3v_4$.
Since $v_1,v_2,v_3$ have no neighbors outside this path,
all other vertices of $T$ can keep their images.

We now prove the first assertion by induction on $n$.
Since $H$ has a vertex of degree at most two, $G$ has a vertex
of degree at least $n-3$. This gives the required path for
$n=4,5$.

For $n=6$, suppose that $G$ has no four-vertex path.
Every connected graph with no such path is a star or a triangle.
Since $G$ has a vertex of degree at least three, one of its
components is a star with at least three leaves.
Choose three of these leaves and a fourth vertex distinct from
them and the center. If this fourth vertex belongs to the star,
it is another leaf; otherwise it lies in a different component.
In either case, the four chosen vertices are pairwise nonadjacent
in $G$, so they form a $K_4$ in $H$, a contradiction.

Finally, let $n\ge7$. Delete any host vertex and apply induction
to obtain a path on $n-3$ vertices in the remaining graph.
Viewed in the original host, this embedding leaves three vertices
unused. Since $n-3\ge4$, we can apply the second assertion to
the last four vertices of the path, adding one new terminal
vertex. The resulting path has $n-2$ vertices.
\end{proof}

\begin{lemma}[A matching surviving one edge deletion]\label{lem:matching}
Let $H$ be a $2$-degenerate graph on $2k+2$ vertices, where $k\ge2$.
Deleting at most one edge from $\overline H$ leaves a graph
containing a matching of size $k$.
\end{lemma}

\begin{proof}
Let $R$ be the graph obtained after the deletion.
Suppose, for a contradiction, that a maximum matching $M$ in $R$
has size at most $k-1$. Let $U$ be the set of vertices not covered
by $M$. Then
\[
 |U|=2k+2-2|M|\ge4.
\]
No two vertices of $U$ are adjacent in $R$, since their edge
could be added to $M$. Hence $H[U]$ is complete except possibly
for the single edge deleted from $\overline H$.

If $|U|\ge5$, any five vertices of $U$ induce in $H$ a graph
of minimum degree at least three, a contradiction.
Thus $|U|=4$. Moreover, an edge must have been deleted, and both
its endpoints lie in $U$; otherwise $H[U]$ would be a $K_4$.
Consequently, $H[U]$ is $K_4$ minus one edge, and $R$ agrees with
$\overline H$ on every pair not contained in $U$.

Since $|U|=4$, we have $|M|=k-1\ge1$. Choose an edge $xy\in M$.
Each of $x$ and $y$ has an $R$-neighbor in $U$. Indeed, if $x$
had none, it would be adjacent in $H$ to all four vertices of $U$.
Then $H[U\cup\{x\}]$ would be $K_5$ minus one edge, contradicting
$2$-degeneracy. The same argument applies to $y$.

If $x$ and $y$ had distinct $R$-neighbors $u,v\in U$, respectively,
we could replace $xy$ in $M$ by $xu$ and $yv$, obtaining a larger
matching. Thus their neighbor sets in $U$ must both consist of
the same single vertex $w$.

Set $W=U\setminus\{w\}$. Both $x$ and $y$ are adjacent in $H$
to every vertex of $W$. Also, $H[W]$ is a triangle with at most
one edge missing, so each vertex of $W$ has a neighbor in $W$.
Therefore, in $H[W\cup\{x,y\}]$, each vertex of $W$ has degree
at least three, while $x$ and $y$ each have the three vertices
of $W$ as neighbors. This contradicts $2$-degeneracy.
\end{proof}

\section{Embedding spiders with a prescribed center}

A \emph{rooted spider} is a tree consisting of $k\ge2$ paths
that share one endpoint, called the \emph{center}, and are otherwise
vertex-disjoint. These paths are its \emph{arms}; their lengths
$a_1,\ldots,a_k$ are positive numbers of edges. Thus the spider
has $1+\sum_{i=1}^k a_i$ vertices.
This definition includes stars and paths rooted at an internal vertex.

\begin{lemma}[Prescribed center]\label{lem:spider}
Let $H$ be a $2$-degenerate graph on $n\ge5$ vertices, and let
$p\in V(H)$ satisfy $d_H(p)\le2$.
Then $\overline H$ contains every rooted spider of order $n-2$,
with its center embedded at $p$.
\end{lemma}

\begin{proof}
Write $G=\overline H$, and let $S$ be the spider to be embedded.
We induct on $n$. For $n=5$, the spider has two arms of length one,
and the assertion follows from $d_G(p)\ge2$.

\smallskip\noindent
\textbf{Reducing long arms.}
Suppose an arm has length at least four.
Delete its terminal leaf from $S$ and any host vertex other than $p$.
The remaining forbidden graph is $2$-degenerate, and the degree
of $p$ in it is still at most two. Induction embeds the smaller
spider with its center at $p$.

In the original host, this embedding leaves three vertices unused.
Apply the terminal repair in Lemma~\ref{lem:paths} to restore the
deleted leaf. Only the three vertices immediately preceding that
leaf may change their images. Since the arm has length at least
four, none of these vertices is the center.
Thus $p$ remains fixed, and we may henceforth assume that
every arm has length at most three.

\smallskip\noindent
\textbf{Partitioning a path into arms.}
Call the vertices in $N_H(p)$ \emph{bad}, and all other vertices
of $G-p$ \emph{good}. There are at most two bad vertices.
A vertex is good precisely when it can be joined to $p$ in $G$.

By Lemma~\ref{lem:paths}, $G-p$ contains a path $P$ on $n-3$
vertices, leaving two vertices of $G-p$ unused.
Since $\sum_{i=1}^k a_i=n-3$, we can partition $P$ into consecutive
blocks of sizes $a_1,\ldots,a_k$.
A block of size $a_i$ with a good endpoint gives an arm of length
$a_i$: join $p$ to that endpoint and then follow the block.
The order of the block sizes along $P$ can be chosen freely.

A singleton block at a bad vertex can be replaced by an unused
good vertex. There are enough such vertices: if $b$ of the two
unused vertices are bad, there are $2-b$ unused good vertices
and at most $2-b$ bad vertices on $P$.
This handles all singleton blocks, and in particular all stars.

A block of size at least two fails to give an arm only when
both endpoints are bad. Thus at most one such block can fail.
If none does, replacing the bad singleton blocks completes
the embedding.

\smallskip\noindent
\textbf{Blocks of different sizes.}
Suppose a block fails and another block has a different size.
Choose a nearest block of a different size. All blocks between
them have the same size as the failing block.
Reorder the block sizes on that side so that the different-sized
block becomes adjacent to the failing one, without changing
the vertices occupied by the failing block.

Now interchange the sizes of these two adjacent blocks and
repartition their combined segment of $P$. The two bad vertices
either lie in different blocks, or one becomes internal to a block.
Indeed, the two sizes are distinct members of $\{1,2,3\}$,
and the bad vertices were the endpoints of the original failing
block. Consequently, neither new block of size at least two
has two bad endpoints. Any bad singleton can again be replaced
by an unused good vertex.

It remains to consider spiders whose arms all have length two
or all have length three.

\smallskip\noindent
\textbf{All arms have length two.}
Here $n-2=1+2k$, so $G-p$ has $2k+2$ vertices.
If there are two bad vertices and they are adjacent in $G$,
delete their edge. Lemma~\ref{lem:matching}, applied to $H-p$,
gives a matching of size $k$ in the resulting graph.
Each matching edge has a good endpoint, since an edge joining
two bad vertices has been excluded. Joining a good endpoint
of each edge to $p$ gives the required arms.

\smallskip\noindent
\textbf{All arms have length three.}
We need only consider a failing block, whose vertices occur
as $a,b,c$ along $P$. Its endpoints $a,c$ are the two bad vertices,
so every other vertex of $G-p$ is good.
Let $u,v$ be the two unused vertices.

If $ac\in E(G)$, replace the block by the path $b,a,c$,
which starts at a good vertex.
Otherwise $ac\in E(H)$.
If $u$ or $v$, say $u$, has a $G$-neighbor $w$ in $\{a,b,c\}$,
choose a neighbor $w'$ of $w$ along the original block.
Then $u,w,w'$ is a replacement path starting at a good vertex.
In both cases all other blocks remain unchanged.

We may therefore assume that $u,v$ are both adjacent in $H$
to all of $a,b,c$. Since $ac\in E(H)$, we must have $uv\in E(G)$;
otherwise $\{a,c,u,v\}$ would form a $K_4$ in $H$.

There is another block because $k\ge2$.
Reverse $P$ if necessary so that a neighboring block follows
$a,b,c$, and denote its vertices by $d,e,f$ in path order.
If $f$ has a $G$-neighbor in $\{u,v\}$, replace the two blocks
by
\[
 (e,d,c)
 \quad\text{and}\quad
 (f,u,v)\ \text{or}\ (f,v,u),
\]
choosing the second path according to that neighbor.
Both paths start at good vertices.

Otherwise $fu,fv\in E(H)$. We claim that $fa\in E(G)$.
If not, the edges
\[
 ac,\quad au,\quad av,\quad cu,\quad cv,\quad
 fu,\quad fv,\quad fa
\]
give a subgraph of $H$ on $\{a,c,u,v,f\}$ with minimum degree
at least three, a contradiction.
We can therefore replace the two blocks by
\[
 (f,a,b)
 \quad\text{and}\quad
 (e,d,c).
\]
Again, both paths start at good vertices. In each construction
the replacement paths are disjoint and avoid all unchanged
blocks. Joining their starting vertices to $p$ completes
the embedding.
\end{proof}

\section{Proofs of the main theorems}

\begin{proof}[Proof of Theorem~\ref{thm:universal}]
We induct on $n$. Let $G=\overline H$, and let $T$ be a tree
of order $n-2$.

For $n=4$, the tree $T$ consists of a single edge.
Since $H$ has a vertex of degree at most two, $G$ has an edge,
as required.

Now let $n\ge5$, and choose $p\in V(H)$ with $d_H(p)\le2$.
If $T$ has at most one vertex of degree at least three,
then $T$ is a spider: its center is that vertex if it exists,
and otherwise we choose any internal vertex of the path $T$.
Lemma~\ref{lem:spider} embeds $T$ into $G$ with its center at $p$.

Suppose, therefore, that $T$ has at least two vertices of degree
at least three, which we call branching vertices.
Root $T$ at one branching vertex and choose a branching vertex
$z$ farthest from the root. Let $w$ be the parent of $z$.
Deleting the edge $zw$ splits $T$ into two components:
let $S$ contain $z$, and let $C$ contain the root.

No proper descendant of $z$ is a branching vertex.
Consequently, $S$ is a spider centered at $z$, with
$d_T(z)-1\ge2$ arms. Write $s=|V(S)|$, so $s\ge3$.
The root has degree at least two in $C$, since deleting $zw$
removes at most one of its incident edges. Thus $|V(C)|\ge3$.

We first embed $C$, reserving $p$ for the center of $S$.
To ensure that the connecting edge $zw$ can later be added,
we embed all of $C$ at vertices adjacent to $p$ in $G$.
There are enough such vertices: writing
$N_H[p]=N_H(p)\cup\{p\}$, we have
\[
 |V(H)\setminus N_H[p]|
 \ge n-3
 \ge n-s
 =|V(C)|+2.
\]
Choose a set
\[
 A\subseteq V(H)\setminus N_H[p],
 \qquad |A|=|V(C)|+2.
\]
The graph $H[A]$ is $2$-degenerate and has fewer than $n$ vertices.
By induction, $C$ has an embedding $\varphi$ into
$\overline{H[A]}$.

Let $B$ be the set of all host vertices not used by this embedding:
\[
 B=V(H)\setminus\varphi(V(C)).
\]
Then
\[
 |B|=n-|V(C)|=s+2.
\]
Also, $p\in B$, since $p\notin A$.
Notice that $B$ includes both the vertices outside $A$ and
the two vertices of $A$ left unused by $\varphi$.

The graph $H[B]$ is $2$-degenerate, and $d_{H[B]}(p)\le2$.
Lemma~\ref{lem:spider} therefore embeds $S$ into
$\overline{H[B]}$ with $z$ at $p$.
The embeddings of $C$ and $S$ use disjoint vertex sets.
Finally, $\varphi(w)\in A\subseteq V(H)\setminus N_H[p]$,
so $p\varphi(w)\in E(G)$. Adding this edge realizes $zw$
and completes the embedding of $T$.
\end{proof}

\begin{proof}[Proof of Theorem~\ref{thm:main}]
Embed the first tree arbitrarily into $K_n$.
When the second tree is revealed, apply Lemma~\ref{lem:pair}
to embed it edge-disjointly from the first, keeping the first
embedding fixed and making their union $H$ $2$-degenerate.
This construction requires no knowledge of the third tree.

When the third tree is revealed, the unused edges form
$\overline H$. By Theorem~\ref{thm:universal}, this graph contains
a copy of the third tree, which has order $n-2$.
Embedding it there completes the packing without changing
either earlier embedding.
\end{proof}

\section*{Declaration on the use of artificial intelligence}

Generative artificial intelligence tools were used during the preparation
of this manuscript to improve the exposition and to assist in developing
and drafting parts of the proof. In particular, several auxiliary
arguments in the proof of Theorem~\ref{thm:main} were initially developed
with the assistance of GPT-5.6. All mathematical statements and proofs
presented here have been checked and revised by the author, who assumes
full responsibility for their correctness.

\end{document}